\documentstyle[12pt]{article}
\begin{document}

\begin{center}
{\Large Exact Random
Walk Distributions using Noncommutative Geometry}
\end{center}

\vspace{1cm}

\begin{center}
Jean Bellissard\dag , Carlos J Camacho\ddag , Armelle 
Barelli\dag , and Francisco Claro\ddag
\end{center}

\vspace{1cm}

\noindent \dag\ Laboratoire de Physique Quantique, UMR5626 associ\'ee au CNRS,
IRSAMC, Universit\'e Paul Sabatier, 118, route de Narbonne, F-31062 Toulouse
Cedex 4, France

\vspace{3mm}

\noindent \ddag\ Facultad de F\'{\i}sica, P. Universidad Cat\'olica de
Chile, Casilla 306, Santiago 22, Chile

\vspace{5mm}

\begin{abstract}
Using the results obtained by the non commutative geometry techniques applied to
the Harper equation, we derive the areas distribution of random walks of length
$ N $ on a two-dimensional square lattice for large $ N $, taking into account 
finite size contributions.
\end{abstract}

\noindent {\bf Pacs Numbers : 05.45.+b, 72.15.Qm, 72.10.Bg}

\vspace{1cm}

Let us consider on a square lattice all closed paths of length $ N $ starting at
the origin. For such a path $ \Gamma $, let $ A(\Gamma ) $ be its algebraic
area. As $ N\mapsto\infty $, the average size of such a path increases as $
\sqrt{N} $ so that $ A(\Gamma )\simeq N $ and the renormalized area will be $
a=A/N $. We want to compute the probability distribution $ {\cal P}(A,N) $ of
the areas at large but finite $ N $. 

\noindent In the limit $ N\mapsto\infty $, the distribution was computed first
by L\'evy using brownian paths \cite{Lev}. We will give a method based upon the
Harper model, allowing to compute the finite size corrections in a systematic
way. The Harper model was designed in 1955 \cite{Har} as the simplest
non-trivial one describing the motion of an electron sitting on a
two-dimensional square lattice and submitted to a uniform magnetic field. Let $
\phi $ be the magnetic flux through the unit cell and let $ \phi _0=h/e $ be the
flux quantum. We set $ \gamma =2\pi\phi /\phi _0 $. Then given $ m=(m_1,m_2)\in
{\bf Z}^2 $, we denote by $ W(m) $ the corresponding magnetic translations
\cite{Zak}. They satisfy the Weyl commutation rules
\begin{equation}
W(m)W(m')=W(m+m')\exp(i{\gamma\over 2}m\wedge m')\,,
\end{equation}
where $ m'\wedge m=m_1'm_2-m_2'm_1\in {\bf Z} $. Note that $ \gamma $ plays a
r\^ole similar to the Planck constant in the canonical commutation relation. 

\noindent Harper's model is given by the following Hamiltonian : 
\begin{equation} 
H=\sum _{\vert a\vert =1}W(a)\mbox{ , }
\end{equation}
where $ \vert a\vert =\vert a_1\vert +\vert a_2\vert $ if $ a=(a_1,a_2)\in {\bf
Z}^2 $. In addition, one defines the trace per unit area as the unique linear
map $ {\cal T} $ on the algebra generated by the $ W(m) $'s such that  
\begin{equation} 
{\cal T}(W(m))=\delta _{m,0}\mbox{ . } 
\end{equation}
Then from (1), (2) and (3), we get :
$$ {\cal T}(H^N)= \sum_{\Gamma\hspace{2mm} :\hspace{2mm}\hbox{closed paths of 
length ${\displaystyle N}$}}
{\rm e}^{i\gamma A(\Gamma )/2}\mbox{ , } $$ 
where the sum is taken on the set of closed paths starting at the origin of
length $ N $. Note that $ N $ should be even to get a non zero sum. Let $ \Omega
_N $ be the number of such closed paths, we then get :
\begin{equation}
\sum_{A=-A_{\rm max}}^{A_{\rm max}}{\cal P}_N (A/N)\exp(ixA/N)=\Omega_N
^{-1}\sum_{\Gamma }\exp(ixA(\Gamma )/N)=
\Omega_N^{-1}{\cal T}(H^N)\vert _{\gamma =x/N}\mbox{ . }
\end{equation}
From this relation we obtain :
\begin{eqnarray}
\Omega_N={\cal T}(H^N)\vert _{\gamma =0}=\int\frac{dk_1dk_2}
{4\pi ^2}\left(
2\cos k_1+2\cos k_2\right)^N\nonumber\\
=\frac{4^{N+1}}{2\pi N}\left(1+{\rm O}
(1/N^2)\right)\hspace{2mm}\mbox{as}\hspace{2mm}N\mapsto\infty
\mbox{ . } 
\end{eqnarray}
Moreover as $ N\mapsto\infty $, for a given value of $ x $, $ \gamma =x/N $
tends to zero, so that we can use a semiclassical argument to compute $ {\cal T}(H^N)
$. It has been shown that the spectrum of $ H $ is made of Landau sublevels
\cite{RaBe} :
\begin{equation}
 E_{\ell }^{\pm }(\gamma)=\pm\left(4-\gamma(2\ell +1)+{\gamma^2\over 16}[1+(2
\ell +1)^2]-{\rm O}
(\gamma^3)\right)\hspace{4mm}\ell =0,1,\ldots ,{\rm O}(1/\gamma )\mbox{ , }
\end{equation}
each with multiplicity per unit area $
g_{\ell }^{\pm }=\gamma /2\pi $. So,
\begin{equation}
{\cal T}(H^N) = \sum_{\pm }\sum_{\ell } \left(
E_{\ell }^{\pm }(\gamma)\right)^N
[\gamma/(2\pi)]\,.
\end{equation}
This gives :
\begin{equation}
{\cal T}(H^N)= {4^{N+1}\over 2\pi N}{x/4\over
\sinh (x/4)}\left[1-\frac{1}{2N}\frac{(x/4)^2}
{\sinh ^2 (x/4)}+{\rm O} (1/N^2)\right]\,.
\end{equation}
Using (4), we obtain the probability distribution :
\begin{equation}
{\cal P}_N(a)={\pi\over \hbox{cosh}^2(2\pi a)}+
{\rm O}(1/N)\,.
\end{equation}
We numerically computed $ {\cal P}_N(a) $ from formula (8) and compared the
result with exact numerical calculations (Figure 1).

\let\picnaturalsize=N
\def\picsize{4.0in}
\def\picfilename{fig1.ps}
\ifx\nopictures Y\else{\ifx\epsfloaded Y\else\input epsf \fi
\let\epsfloaded=Y
\centerline{\ifx\picnaturalsize N\epsfxsize \picsize\fi
\epsfbox{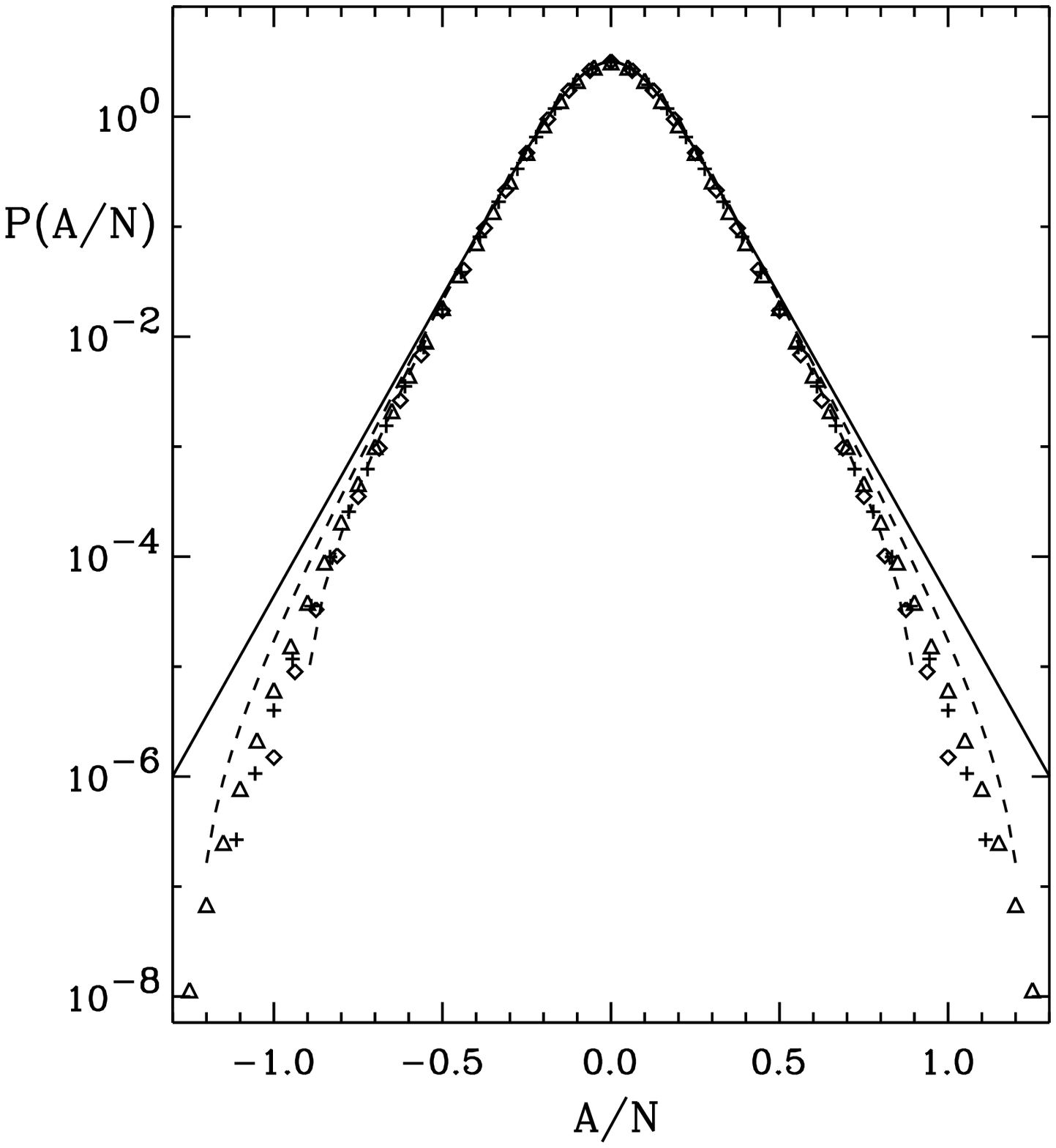}}}\fi

FIG. 1. Scaling function for the probability of having a loop of area A for a
random walk of N steps. Diamond, plus and triangle symbols correspond 
to the finite size data for N = 16, 18 and 20, respectively.
Solid line corresponds to the universal function
(9), as $ N\mapsto\infty $; whereas dashed lines include the $ 1/N $ 
correction term for N=20 and 40. 

\vspace{1cm}

\noindent {\bf Acknowledgements} : This work has been partially supported by the
ECOS/CONYCIT project C94E07 and FONDECYT No. 3940016 (Chile).

\end{document}